\documentclass[12 pt]{article}
\usepackage{amsmath,amssymb}
\usepackage[english]{babel}
\usepackage{lipsum}
\usepackage{hyperref}
\date{}
\begin{document}
\title{\textbf{A SYNTHETIC OVERVIEW ON SOME KNOWN CHARACTERIZATIONS OF WOODIN CARDINALS}}
\author{\sc Gabriele Gull\`a}
\maketitle
\begin{abstract} 
This brief survey comes from the slides of a seminar I gave to philosophy of mathematics students. I will present some different characterizations of Woodin cardinals, including the one obtained by Ernest Schimmerling in \cite{S}.\\
I will try to give to this paper the most self-contained possible structure, also by showing explicitly just the proofs usually leaved to the reader, and giving exact references for the others.
\end{abstract}
\textbf{Key Words}: Large cardinals, Mostowski collapse, Skolem hull, Extender, Woodin cardinals, weakly hyper-Woodin and hyper-Woodin cardinals, Shelah cardinals.
\let\thefootnote\relax\footnotetext{2020 Mathematical Subject Classification: 03-01, 03E55}
\section{Basic notions}
Let us start by recalling some basic definitions:\\
\\\textbf{Definition 1.1}\\
A function $j$ between two structures $N$ and $M$ is an \textit{elementary embedding} ($\mathcal{EE}$) if for every first order formula $\phi(\underline{x})$ and an element $\underline{a}$ in $N$, it holds that 
$$N\models \phi(\underline{a})\Leftrightarrow M\models \phi(\underline{j(a)})$$
If $N$ and $M$ are transitive classes, we say that $\alpha$ is the \textit{critical point} of $j$ ($\alpha=crit(j)$) if it is the smallest ordinal for whom $j$ is not the identity.\\
\\\textbf{Definition 1.2}\\
Let $\kappa$ be a cardinal; with the symbol $H_{\kappa}$ we identify the collection of the \textit{sets hereditarily of cardinality less $\kappa$}, meaning that they are of cardinality less than $\kappa$, and all members of their transitive closure are of cardinality less than $\kappa$.\\
\\\textbf{Definition 1.3} (Woodin cardinal)\\
Let $\kappa<\lambda$ be two cardinals and let $X$ be a set. We say that\\ $\kappa$ is $(\lambda-X)-$strong if and only if there is a \textit{transitive class} $M$\\ (i.e. $x\in M\Rightarrow x\subseteq M$) and an elementary embedding $j:V\rightarrow M$ such that
$$\kappa=crit(j)$$
$$j(\kappa)\ge \lambda$$
$$j(X)\cap H_{\lambda}=X\cap H_{\lambda}$$
If $\kappa<\gamma$ then $\kappa$ is $(<\gamma-X)-$strong if and only if $\kappa$ is $(\lambda-X)-$strong for every $\lambda<\gamma$.\\
A cardinal $\delta$ is a \textit{Woodin cardinal} if and only if it is strongly inaccessible ($2^{\lambda}<\delta$ for every $\lambda<\delta$) and for every $X\subseteq H_{\delta}$ there is a $\kappa<\delta$ which is $(<\delta-X)-$strong.\\
\\\hypertarget{+}{\textbf{Definition 1.4}} (Extender)\\
Let $j:V\rightarrow M$ be an $\mathcal{EE}$ with  $crit(j)=\kappa\le \lambda\le j(\kappa)$; for every finite subset $A\subset \lambda$ we define a measure $E_A$ on $[\kappa]^{<\omega}$ as follows
$$X\in E_A\Leftrightarrow A\in j(X)$$ 
so we define an \textit{extender} $\mathbb{E}$ as the collection of measures
$$\left\{E_A: A\in [\lambda]^{<\omega}\right\}$$
\textbf{Remark 1.1}\\
It is obvious that this extender depends on $\kappa, \lambda$ and $j$: we avoid to underline this dependence in our notation.\\
\\Now we define the quotient
$$Ult_{\mathbb{E}}=\left\{[A, f]_{\mathbb{E}}: A\in [\lambda]^{<\omega}, f:[\kappa]^{|A|}\rightarrow V\right\}$$
where we identify $(A, f)$ with $(B, g)$ if and only if $$\left\{t\in [\kappa]^{|A\cup B|}: f(\pi_{A\cup B, A}(t))=g(\pi_{A\cup B; B}(t))\right\}\in E_{A\cup B}$$ (where $\pi_{C, D}:[\lambda]^{|D|}\rightarrow [\lambda]^{|C|}$, with $D\supseteq C$, maps $\left\{x_1...x_n\right\}$ in $\left\{x_{i_1}....x_{i_m}\right\}$).\\
Finally we define $j_{\mathbb{E}}:V\rightarrow Ult_{\mathbb{E}}$ as the function which maps $Y$ in $[\emptyset, c_Y]$, where $c_Y$ is the function constantly equal to $Y$.\\
\\\textbf{Theorem 1.1}\\The following are equivalents:\\
1) $\kappa$ is a Woodin cardinal.\\
2) $\forall f\in  ^{\kappa}\kappa\enspace \exists \alpha<\kappa|\enspace f''\alpha\subseteq \alpha\enspace\wedge\enspace\exists\enspace j:V\rightarrow M|\enspace crit(j)=\alpha\enspace \wedge\enspace V_{j(f)(\alpha)}\subseteq M$.\\
3) $\forall A\subseteq V_{\kappa}$
$$\left\{\alpha<\kappa|\enspace \alpha\enspace (\gamma-A)-\mbox{strong}\enspace \forall \gamma<\kappa\right\}$$
is \textit{stationary} in $\kappa$ (it intersects every $C\subseteq\kappa$ such that $\sup{(C)}=\kappa$ and $\forall\enspace\mbox{limit}\enspace\gamma<\kappa\enspace (\sup(C\cap \gamma)=\gamma\Rightarrow \gamma\in C)$).\\
4) $F=\left\{X\subseteq \kappa\enspace|\enspace \kappa - X\enspace \mbox{is not Woodin in}\enspace \kappa\right\}$ is a proper filter over $\kappa$ (the \textit{Woodin filter}).\\
5) $\forall f\in\enspace  ^{\kappa}\kappa\enspace \exists \alpha<\kappa|\enspace f''\alpha\subseteq \alpha\enspace\wedge\enspace\exists\enspace \mbox{an extender}\enspace \mathbb{E}\in V_{\kappa}\enspace |$ $$crit(j_{\mathbb{E}})=\alpha$$
$$j_{\mathbb{E}}(f)(\alpha)=f(\alpha)$$ 
$$V_{j_{\mathbb{E}}(f)(\alpha)}\subseteq Ult_{\mathbb{E}}$$
\\\textbf{Some words on the proof and references for a complete one:}\\
Obviously \textbf{(3)$\Rightarrow$(1)}, because, if the required $\alpha$ did not exist, then the set defined in (3) would be empty and so, by definition, not stationary.\\
\textbf{(5)$\Leftrightarrow$(2)} thanks to the functions $j\rightarrow E_j\wedge E\rightarrow j_E$ in \hyperlink{+}{Definition 1.4}.\\
For \textbf{(4)$\Leftrightarrow$(2)}: if $F$ is a proper filter then $\emptyset\notin F$ and so $\kappa$ is a Woodin cardinal.\\
\textit{Viceversa}, let $\kappa$ be a Woodin cardinal, then obviously $\emptyset$ is in $F$; if $A, B\in F$ it is not possible to find a required $\alpha$ in $\kappa-(A\cap B)$ so $A\cap B$ is in $F$ too, and similarly if $B\supset A\in F$ such an $\alpha$ will not be in $\kappa-B$, and so $B\in F$ too.\newpage\noindent
The remaining part of the proof can be found in \cite{J} or \cite{Ka}; in particular for \textbf{(2)$\Rightarrow$(3)} see proposition 26.13 and the first part of Theorem 26.14 in \cite{Ka}; the proof of this last Theorem shows also \textbf{(1)$\Rightarrow$(5)}.\\ 
\begin{flushright}$\mathcal{a}$\end{flushright}
\section{The Schimmerling characterization}
We open this section with the preliminary results needed in order to present the characterization that Ernest Schimmerling gives of Woodin cardinals in \cite{S}.\\
\\\textbf{Lemma 2.1} (Mostowski Collapse, MC)\\
Let $E$ be a binary relation on a class $X$ such that:\\
1) $E$ is \textit{set-like}: $\left\{y\enspace|\enspace yEx\right\}$ is a set for every $x\in X$;\\
2) $E$ is \textit{well-founded}: every non empty subset of $X$ contains an $E$-minimal element;\\
3) the structure $\left\langle X, E\right\rangle$ is \textit{extensional}:\\ $\forall x, y\in X\enspace [(zEx\Leftrightarrow zEy\enspace \forall z\in X)\Rightarrow  x=y]$;\\
then there are a unique isomorphism $\pi$ and a unique transitive class  $M$ such that $$\left\langle X,E \right\rangle \simeq^{\pi} \left\langle M, \in \right\rangle$$
e $\pi(x)= \left\{\pi(y)\enspace|\enspace y\in X\wedge  yEx\right\}$ (where obviously $\pi$ is defined by recurrence thanks to the well-foundedness of $E$: the ``0 step'' is the minimal element).\\
\\For the proof see, for example, Theorem 6.15 in \cite{J} or Lemma I.9.35 in \cite{Ku}.\\\begin{flushright}$\mathcal{a}$\end{flushright}
\textbf{Definition 2.1} (Skolem Hull, SH)\\
Given a first order formula, we call \textit{skolemization} the replacement process of $\exists $-quantified variables with terms of the type $f(\underline{x})$.\\
The (new) symbol ``$f$'' identifies a \textit{Skolem function}.\\
A theory that, for every formula with free variables $\underline{x}, y$ has a Skolem function is called \textit{Skolem Theory}.\\
Given a model $\mathfrak{M}$ of a Skolem theory and a set $X$, the smallest substructure containing $X$ is called \textit{Skolem hull} of X.\\
\\\\\textbf{Definition 2.2} (Schimmerling)\\
Let $M$ be a transitive class and let $\pi:M\rightarrow H_{\theta}$ be an $\mathcal{EE}$ with $\kappa=crit(\pi)$ and let $\lambda<\pi(\kappa)$; finally let $j:V\rightarrow N$ be another $\mathcal{EE}$: then \textit{$j$ certifies $\pi$ up to $\lambda$} if and only if
$$\kappa=crit(j)$$
$$j(\kappa)\ge \lambda$$
$$j(A)\cap H_{\lambda}=\pi(A)\cap H_{\lambda}$$
for every $A\in \mathcal{P}(H_{\lambda})\cap M$.\\
We say that $\pi$ is \textit{certified} if and only if for every $\lambda<\pi(\kappa)$ there is an $\mathcal{EE}$ $j:V\rightarrow N$ which certifies $\pi$ up to $\lambda$.\\
\\\hypertarget{*}{\textbf{Proposition 2.1}} (Schimmerling)\\
Let $M$ be a transitive set, $\pi:M\rightarrow H_{\theta}$ a non trivial $\mathcal{EE}$, $\kappa=crit(\pi)$ and $\lambda<\pi(\kappa)$.\\
Let us suppose that $j:V\rightarrow N$ certifies $\pi$ up to $\lambda$ and let $S$ be an element of the image of $\pi$. Then $\kappa$ is $(\lambda-S)-$strong (witnessed by $j$).\\
\\\textbf{Proof}:\\
By definition of certified $\mathcal{EE}$, in order to show that $\kappa$ is $(\lambda-S)-$strong it suffices to prove that $j(S)\cap H_{\lambda}=S\cap H_{\lambda}$: about this we observe that \textbf{(1)} $S\cap H_{\kappa}\in M$, as $\kappa$ is the first ordinal moved by $\pi$, and obviously \textbf{(2)} $\pi(S\cap H_{\kappa})=S\cap H_{\pi(\kappa)}$.\\
Now,
$$\textbf{(a)}\quad j(S)\cap H_{\lambda}=j(S\cap H_{\kappa})\cap H_{\lambda}$$
since $j(\kappa)\ge \lambda$; moreover, by the hypothesis
$$\textbf{(b)}\quad j(S\cap H_{\kappa})\cap H_{\lambda}=\pi(S\cap H_{\kappa})\cap H_{\lambda}$$
because, if $S\cap H_{\kappa}=A$, then $A\subseteq H_{\kappa}$, so $A\in \mathcal{P}(H_{\kappa})$, and if $A\in M$ it follows that $A\in \mathcal{P}(H_{\kappa})\cap M$.\\
From \textbf{(2)} it follows that
$$\textbf{(c)}\quad \pi(S\cap H_{\kappa})\cap H_{\lambda}=S\cap H_{\lambda}$$
and so we obtain the equivalence between the first member of \textbf{(a)} and the second of \textbf{(c)}:
$$j(S)\cap H_{\lambda}=S\cap H_{\lambda}$$
\begin{flushright}$\mathcal{a}$\end{flushright}
\hypertarget{**}{\textbf{Proposition 2.2}} (Schimmerling)\\
Let $j:V\rightarrow N$ be an $\mathcal{EE}$ with $crit(j)=\kappa$. Let $\theta>\kappa$ be a cardinal and $S\in H_{\theta}$.\\
Let us suppose that $\pi:M\rightarrow j(H_{\theta})$ be the inverse of the MC of the SH of $\kappa\cup \left\{j(S)\right\}$ in $j(H_{\theta})$.\\
Then
$$\textbf{(a)}\quad \kappa=crit(\pi)$$
$$\textbf{(b)}\quad \pi(\kappa)\ge j(\kappa)$$
$$\textbf{(c)}\quad j(A)=\pi(A)\cap j(H_{\kappa})$$
for every $A\in \mathcal{P}(H_{\kappa})\cap M$.\\
\\\textbf{Proof}\\
\textbf{(a)} Let us suppose that $\pi(\kappa)=\kappa$; $\pi(\kappa)\in j(H_{\theta})$ then $\kappa=\pi(\kappa)=j(\gamma)$ but $\kappa=crit(j)$ and so it can not exist a $\gamma<\kappa$ which is moved by $j$. Similarly one shows \textbf{(b)}: if $\pi(\kappa)$ was smaller than $j(\kappa)$ then it would be the image, through $j$, of some $\gamma<\kappa$, which is impossible.\\
\textbf{(c)} $\pi(A)$ e $j(H_{\kappa})$ are in $j(H_{\theta})$, so $\pi(A)\cap j(H_{\kappa})\subseteq j(H_{\theta})$.\\
Clearly $j^{-1}(\pi(A)\cap j(H_{\kappa}))\subseteq H_{\kappa}$, and since $A\in \mathcal{P}(H_{\kappa})$ and $\kappa=crit(j)$ we have that $j^{-1}(\pi(A)\cap j(H_{\kappa}))=H_{\kappa}\cap(\pi(A)\cap j(H_{\kappa}))=A$ (because $\kappa$ is critical point for $\pi$ too).\\By applying $j$ to both members we obtain \textbf{(c)}.\\
\begin{flushright}$\mathcal{a}$\end{flushright}
\hypertarget{***}{\textbf{Proposition 2.3}} (Schimmerling)\\
Let $\delta$ be a Woodin cardinal. Let $\theta>\delta$ be a cardinal and $S\in H_{\theta}$. Let $T$ be the first order theory of $\delta\times \left\{S\right\}$ in $H_{\theta}$ coded as a subset of $H_{\delta}$. Let be $\kappa$ $(<\delta-T)-$strong and $\pi$ the inverse of the MC of the SH of $\kappa\cup S$ in $H_{\theta}$.\\
Then
$$\kappa=crit(\pi)$$
$$\pi(\kappa)\ge \delta$$
Moreover if $\lambda<\delta$ and $j:V\rightarrow N$ makes $\kappa$ $(\lambda-T)-$strong, then $j$ certifies $\pi$ up to $\lambda$.\\
\begin{flushright}$\mathcal{a}$\end{flushright}
The proof (see Proposition 3, Lemma 3.1 and Lemma 3.2 in \cite{S}) shows that $X\cap \delta\subseteq \kappa$ (the other inclusion is trivial because $\kappa$ is included in both $X$ and $\delta$) and, using \hyperlink{**}{Proposition 2.2} and $(\lambda-T)-$strongness, that $j$ certifies $\pi$ up to $\lambda$.\\
\\It implies the following characterization\\
\\\textbf{Theorem 2.1} (Schimmerling)\\
Let $\delta$ be an inaccessible cardinal. Then the following are equivalent:\\
i) $\delta$ is a Woodin cardinal;\\
ii) for every $S\in \delta$ there is a $\kappa<\delta$ such that for every cardinal $\theta>\delta$, if $\pi$ the inverse of the MC of the SH of $\kappa\cap \left\{\delta, S\right\}$ in $H_{\theta}$, then $\pi(\kappa)=\delta$ and $\pi$ is certified.\\
\\\textbf{Proof}:\\
If (ii) holds then hypothesis of \hyperlink{*}{Proposition 2.1} are verified: there is a $j$ which witnesses that $\kappa$ is $(\lambda-S)-$strong and moreover that it is ($\delta-S)-$strong, so $\delta$ is a Woodin cardinal.\\
\textit{Viceversa}, if $\delta$ is a Woodin cardinal and $\pi$ is the inverse of the MC of the SH of $\kappa\cup\left\{\delta, S\right\}$, then we can apply  \hyperlink{***}{Proposition 2.3}: the first part says that $\pi(\kappa)=\delta$, and the second that $\pi$ is certified (being $\delta$ ($\lambda-S$)-strong for every $\lambda<\delta$).\\\begin{flushright}$\mathcal{a}$\end{flushright}
The interested reader who approaches this topics for the first time can complete the overview about ``Woodin-realated'' cardinals (and more) with the last part of \cite{S}, about two new notions of large cardinals connected to the Woodin's one:\\
\\\\\textbf{Definition 2.3}\\
A cardinal $\delta$ is a \textit{weakly hyper-Woodin cardinal} if and only if for every set $S$ there is an ultrafilter $U$ over $\delta$ such that 
$$\left\{\kappa<\delta|\enspace \kappa\enspace\mbox{is}\enspace (<\delta-S)-\mbox{strong}\right\}\in U$$
A cardinal $\delta$ is \textit{hyper-Woodin} if $U$ does not depend on $S$, meaning that  $\delta$ is a Woodin cardinal and $U$ extends the Woodin filter.\\
\\Schimmerling, thanks also to an observation due to Cummings, showed that if $\delta$ is a Shelah cardinal then it is weakly hyper-Woodin.\\ It implies that the increasing hierarchy of this ``Woodin-related'' cardinals is the following:\\
\\\centerline{measurable Woodin $\nearrow$ weakly hyper-Woodin $\nearrow$ Shelah $\nearrow$ hyper-Woodin}\\
\\For further details see what follows Theorem 5 in \cite{S}.\\\\

\end{document}